\documentclass[12pt]{amsart}
\usepackage{amsfonts}
\usepackage{amsmath}
\usepackage{amssymb}
\usepackage{a4}
\usepackage{multirow}
\usepackage{hyperref}

\newcommand{\heute}{August 12th, 2012}

\theoremstyle{plain}
\newtheorem{theorem}{Theorem}[section]

\newtheorem{lemma}[theorem]{Lemma}

\theoremstyle{remark}

\newtheorem*{defn}{Definition}

\newcommand{\ignore}[1]{}

\newcommand{\f}[1][p]{\mathbb{F}_{#1}}

\newcommand{\Gro}[1]{Gr\"ob\-ner}

\newcommand{\depth}{\operatorname{depth}}

\newcommand{\Reg}{\operatorname{Reg}}

\DeclareMathOperator{\Bild}{Im}
\DeclareMathOperator{\Res}{Res}

\DeclareMathSymbol\normal{\mathrel}{AMSa}{"43}

\newcommand{\prank}[1][p]{\operatorname{rk}_{#1}}

\newcommand{\gendeg}{\operatorname{gendeg}}
\newcommand{\tr}{\operatorname{tr}}

\begin{document}

\title[Completeness criteria for modular group cohomology]{Comparing
  completeness criteria for modular cohomology rings of finite groups}
\author[S. A. King]{Simon A. King} \thanks{This work was supported by Marie
  Curie grant MTKD-CT-2006-042685 and DFG grant GR 1585/6--1.}
\address{Department of Mathematics and Computer
  Science\\ Friedrich-Schiller-Universit\"at, D-07737 Jena, Germany}
\email{simon.king@uni-jena.de}
\date{\heute}

\begin{abstract}
\noindent
We introduce a criterion for the completeness of ring approximations of
modular cohomology rings of finite non prime power groups, and discuss how
one can benefit from a combination of different completeness criteria in
practical computations.
\end{abstract}

\maketitle
\newlength{\djglength}

\section{Introduction}
\noindent
Let $G$ be a finite group, and let $p$ be a prime dividing $|G|$. We want to
compute a minimal ring presentation of the cohomology ring $H^*(G;\f[p])$. A
\emph{ring approximation} of $H^*(G;\f[p])$ out to degree $n$ consists of a
graded commutative ring $\tau_nH^*(G;\f[p])$ and a graded ring homomorphism
$\alpha_n: \tau_nH^*(G;\f[p])\to H^*(G;\f[p])$ whose restriction on degree $d$
is an isomorphism of vector spaces, for any $d\le n$.

We assume that we are able to compute a ring approximation of $H^*(G;\f[p])$
out to any finite degree. 
Since any modular cohomology ring of a finite group has a finite presentation,
$\alpha_n$ is a graded ring isomorphism, provided that $n$ is large
enough. Let $n_0$ be the smallest number such that $\alpha_{n}$ is an
isomorphism for all $n\ge n_0$.  A \emph{completeness criterion} is an
algorithmic procedure that can decide whether $\alpha_n$ is an isomorphism. To
be precise: 
\begin{enumerate}
\item Effectiveness: There is some number $N$ such that the criterion asserts
  that $\alpha_n$ is an isomorphism, for all $n\ge N$.
\item Correctness: Let $N_0$ be the smallest number such that the criterion
  asserts that $\alpha_{N_0}$ is an isomorphism. Then $N_0\ge n_0$.
\end{enumerate}

For practical computations, it is important that $N_0-n_0$ is small. And of
course, it is also important that the computations, which the completeness
criterion relies on, are not too difficult.

J.~F.~Carlson~\cite{Carlson:Tests} proposed the first completeness criterion.
Even though the effectiveness of his criterion relies on an unproven
conjecture, it was an essential ingredient of the first complete computation
of modular cohomology rings for all 267 groups of order
64~\cite[Appendix]{CarlsonTownsley}.\@

D.~J.~Benson gave a completeness criterion~\cite{Benson:DicksonCompCoho} that
does not rely on a conjecture.  For Benson's criterion, one needs to construct
elements $\zeta_1,...,\zeta_r\in\tau_nH^*(G;\f[p])$ in degrees at least two
whose images under $\alpha_N$ provide a filter regular homogeneous system of
parameters (f.~r.~hsop, for short) of $H^*(G;\f[p])$, for sufficiently large
$N\ge n$.

Benson suggests to construct the f.~r.~hsop using Dickson invariants in the
cohomology rings of $p$-elementary abelian subgroups of $G$. The resulting
degrees of the parameters grow exponentially in the $p$-rank $\prank[p](G)$,
and one has $N_0>\sum(|\zeta_i|-1)$ for this criterion. Moreover, it involves
the computation of the kernels of the multiplications maps associated with
$\zeta_i$ to determine the so-called \emph{filter degree type}.  Both the
large degree and the computation of the filter degree type can be problematic.

With D.~Green~\cite{128gps}, we modified Benson's criterion, so that it can
apply much earlier. On the one hand, it is based on an improved construction
of a filter regular hsop, often yielding fairly small degrees. On the other
hand, it uses an existence result for filter regular parameters in even
smaller degrees over a finite extension field. This criterion was involved in
the first computation of the modular cohomology rings of all groups of order
128~\cite{128gps}, of the Sylow $2$-subgroup of the third Conway
group~\cite{GreenKingEllis:ConwayThree}, and of all but six groups of order
243~\cite{GreenKing:128website}.
However, it is still needed to explicitly construct a f.~r.~hsop and compute
the filter degree type, which in some cases is very difficult. For instance,
the computation of the filter degree type of the mod-$2$ cohomology of
$N_\mathrm{HS}(Z(\mathop{Syl}_2\mathrm{HS}))$, where $\mathrm{HS}$ is the
Higman-Sims group, took several days of computation time on a decent computer.

P.~Symonds~\cite[conluding remark]{Symonds:RegProof} suggests a different
criterion. It relies on constructing homogeneous elements
$\zeta_1,...,\zeta_r\in\tau_nH^*(G;\f[p])$ whose images in
$H^*(G;\f[p])$ generate a sub-algebra over which $H^*(G;\f[p])$ is a finitely
generated module. By abuse of the usual notion, we refer to such elements
as \emph{parameters}, even though they may be algebraically dependent.

The Symonds criterion detects completeness if $n > \sum(|\zeta_i|-1)$ and
$\tau_nH^*(G;\f[p])$ is generated in degree at most $n$ as a module over the
parameters. While it is not needed to compute the filter degree type, Symonds'
criterion requires the explicit construction of parameters. If the parameter
degrees are fairly high then $N_0$ can be much larger than $n_0$, which would
be a problem. But the construction of parameters in small degrees can be a
problem as well.
% By lack of a
%different construction, one would probably still start with a f.-r.~hsop as
%in~\cite{Benson:DicksonCompCoho} or~\cite{128gps}\@, and then simplify it. 
It seems impossible to improve the criterion by exploiting an existence proof
for parameters over finite extension fields as in~\cite{128gps}. However, once
the parameters are constructed, the application of Symonds' criterion is
easy.

In this paper, we suggest a completeness criterion for $H^*(G;\f[p])$ in the
case that $G$ is a finite group that is not of prime power order. It relies on
knowing $H^*(U;\f[p])$ for a subgroup $U<G$ whose index in $G$ is co-prime to
$p$. When using the stable element method (see, \emph{e.g.}, \cite[XII
  \S10]{CartanEilenberg}), knowledge of $H^*(U;\f[p])$ is granted anyway.
% Completeness is then proved by exploiting the Poincar\'e series.

The criterion, that we call the \emph{Hilbert--Poincar\'e criterion}, has two
parts. The first part, namely Lemma~\ref{lem:gen}, is used to test whether
$\tau_nH^*(G;\f[p])$ contains all generators of $H^*(G;\f[p])$, hence, whether
$\alpha_n$ is surjective. It relies on the assumption that $G$ is not of prime
power order and that we know $H^*(U;\f[p])$.

The second part of the criterion, namely Theorem~\ref{thm:rel}, applies to any
finite group (including those of prime power order), but it relies on knowing
surjectivity of $\alpha_n$. It uses the existence of parameters in small
degrees for the cohomology ring of $G$ with coefficients in a finite extension
field. It then relates the Poincar\'e series of $\tau_{n}H^*(G;\f[p])$ with
the degrees of these parameters and with the depth of $H^*(U;\f[p])$ to test
whether $\alpha_n$ is an isomorphism. Lemma~\ref{lem:pars} provides one way to
prove the existence of small parameters.

Hence, the Hilbert--Poincar\'e criterion combines the advantages of the
improved Benson criterion (the existence of small parameters over an extension
field can be used) and the Symonds criterion (the application of the criterion
only involves a relatively easy computation).

Even if a different criterion is used to detect completeness,
Lemma~\ref{lem:gen} is very handy when computing a cohomology ring
approximation with the stable element method. Here, $H^*(G,\f[p])$ is
identified with a graded sub-ring of $H^*(U,\f[p])$, with a subgroup $U<G$
containing a Sylow $p$-subgroup of $G$. For any degree $d$, $H^{(d)}(G,\f[p])$
can be computed by solving a system of linear equations (the \emph{stability
  conditions}) in $H^{(d)}(U,\f[p])$.
% To compute the $n$-th ring
%approximation, one would usually start with the power products of generators
%of $H^*(G,\f[p])$ in degrees less than $n$, which allows to find new algebraic
%relations between these generators. Then, one would compute $H^(n)(G,\f[p])$
%by solving the stability conditions, would compare it with the subspace
%spanned by power products, and may find further generators of
%$H^*(G,\f[p])$. But 
However, if Lemma~\ref{lem:gen} asserts that $\alpha_n$ is surjective, then
solving the stability conditions is not needed, for $d\ge n$.

Section~\ref{sec:surjective} is devoted to proving
Lemma~\ref{lem:gen}. Section~\ref{sec:complete} proves the second part of the
Hilbert--Poincar\'e criterion. Section~\ref{sec:example} illustrates the
benefits of the first part of the Hilbert--Poincar\'e criterion by an example
and suggests a heuristics to make best use of the available completeness
criteria.

\subsection*{Acknowledgment}
We are grateful for interesting e-mail exchange with Peter Symonds. We owe to
him the idea of using the Poincar\'e series for a completeness criterion.

\section{Surjectivity of the ring approximation}
\label{sec:surjective}

Let $p$ be a prime number. Let $G$ be a finite group whose order is divisible
by $p$ but is not a prime power. Let $U<G$ be a proper subgroup such that $p$
does not divide $[U:G]$. The embedding of $U$ in $G$ induces the restriction
map $\Res^G_U: H^*(G;\f[p])\to H^*(U;\f[p])$. On the other hand, we have the
transfer map $\tr^U_G:H^*(U;\f[p])\to H^*(G;\f[p])$.

Let us recall some well known facts available, \emph{e.g.},
in~\cite{Evens:book}. For $x\in H^*(U;\f[p])$ and $y\in H^*(G;\f[p])$ holds
$\tr^U_G(\Res^G_U(y)) = [G:U]\cdot y$ and $\tr^U_G(\Res^G_U(y)\cdot
x)=y\cdot\tr^U_G(x)$. Since $[G:U]$ is invertible in $\f[d]$, it follows that
$\Res^G_U$ is injective.

In the case where $U$ is  Sylow $p$-subgroup of~$G$, \cite[XII
  \S10]{CartanEilenberg} show that the image of $\Res^G_U$ is characterised by
so-called stability conditions, that are associated with double cosets of
${U\setminus G}/U$. The method of proof works equally well in the general
case~\cite[Prop.\@ 3.8.2]{Benson:I}. Hence, if $H^*(U;\f[p])$ is known then
one can compute $H^{(n)}(G;\f[p])$ by solving systems of linear equations in
$H^{(n)}(U;\f[p])$, for any $n$.

\begin{defn}
We consider $H^*(U;\f[p])$ as a $\tau_nH^*(G;\f[p])$-module via $\Res^G_U\circ
\alpha_n$. Let
$$\gendeg_n(G,U) = \min \left\{ d\in \mathbb N: H^*(U;\f[p]) =
\Bild\left(\Res^G_U\circ\alpha_n\right)\cdot H^{\le d}(U;\f[p])\right\}$$
\end{defn}

Note that $\gendeg_n(G,U)=\infty$ if $H^*(U;\f[p])$ is not finitely generated
as a $\tau_nH^*(G;\f[p])$-module. Also note that $\gendeg_{n+1}(G,U)\le \gendeg_n(G,U)$.

\begin{lemma}\label{lem:gen}
If $n$ is sufficiently large then $\gendeg_n(G,U)$ is finite. If
$n\ge\gendeg_n(G,U)$, then $\alpha_n$ is surjective.
\end{lemma}
\begin{proof}
It is well known that $H^*(U;\f[p])$ is finitely generated as a
$H^*(G;\f[p])$-module via restriction. If $n$ is big enough then $\alpha_n$ is
an isomorphism, and the first statement follows.

We prove the second assertion by contradiction. Let $n\ge\gendeg_n(G,U)$, and
assume that there is some $y\in H^*(G;\f[p])\setminus \Bild(\alpha_n)$. Since
$\alpha_n$ is an isomorphism out to degree $n$, we have $|y|>n\ge
\gendeg_n(G,U)$. By definition of $\gendeg_n(G,U)$, there are $x_1,...,x_k\in
H^{\le \gendeg_n(G,U)}(U;\f[p])$ and $y_1,...,y_k\in \tau_nH^*(G;\f[p])$ such
that
$$ \Res^G_U(y) = \sum_{i=1}^k \Res^G_U(\alpha_n(y_i))\cdot x_i.$$ Hence,
$$ [G:U]\cdot y = \tr^U_G(\Res^G_U(y)) = \sum_{i=1}^k \alpha_n(y_i)\cdot
\tr^U_G(x_i).$$ Since $|x_i|\le \gendeg_n(G,U)\le n$ and $\alpha_n$ is an
isomorphism out to degree $n$, $\tr^U_G(x_i)\in\Bild(\alpha_n)$. Since $[G:U]$
is invertible in $\f[p]$, it follows that $y\in\Bild(\alpha_n)$.
\end{proof}

\section{Completeness of the ring approximation}
\label{sec:complete}

Let $R$ be a finitely generated graded commutative $\f[p]$-algebra. If
$X\subset R$ is a set of homogeneous elements, we denote by $\langle
X\rangle\subset R $ the two sided ideal generated by $X$. We denote the
Poincar\'e series of $R$ by $P(R;t)$. Denote $R_L = L\otimes_{\f[p]}R$ for any
extension field $L$ of $\f[p]$, and consider $R\subset R_L$ by slight abuse of
notation.

Let $\mathcal P=\{\zeta_1,...,\zeta_r\}\subset R$, and assume that $\mathcal
P$ is a set of parameters, \emph{i.e.}, $R$ is a finitely generated module
over the subalgebra spanned by $\mathcal P$. We denote the
Krull dimension~\cite{Benson:MSRI} of $R$ by $\dim (R)$.
% We are interested in making the
%degrees as small as possible. So, for any $i=1,...,r$, we can assume that
%$R/\left\langle \mathcal P\setminus\{\zeta_i\}\right\rangle$ has no parameter
%in degree $<|\zeta_i|$ --- otherwise, we could easily reduce the degrees in
%$\mathcal P$.

\begin{lemma}\label{lem:pars}
Let $X\subset \mathcal P$. Assume that $R/\langle X\cup R^{(d)}\rangle$ is a
finite-dimensional $\f[p]$-vector space, for some $d\in\mathbb N$. Then there
exists a finite extension field $K$ of $\f[p]$ so that $R_K$ has a hsop formed
by $X$ together with $\dim(R/\langle X\rangle)$ elements of degree
$d$.
\end{lemma}

\begin{proof}
Let $\tilde R=R/\langle X\rangle$. 
%Since $X$ is subset of a hsop of $R$, we
%have $\dim \tilde R=\dim R - |X|$.
From the assumption follows that $\tilde R/\langle \tilde R^{(d)}\rangle$ is a
finite dimensional $\f[p]$-vectorspace.

Let $L$ be an infinite field that is an algebraic extension field of
$\f[p]$. Any $\f[p]$-basis $B_d$ of $\tilde R^{(d)}$ can also be interpreted
as an $L$-basis of $\tilde R_L^{(d)}$.  Since $\tilde R_L/\langle \tilde
R_L^{(d)}\rangle$ is a finite dimensional $L$-vectorspace and $L$ is infinite,
one version of Noether normalisation allows to conclude the existence of
$L$-linear combinations $p_{i}$ over $B_d$, forming a hsop of $\tilde R_L$,
where $i=1,...,\dim \tilde R_L$.

Now, let $K$ be the smallest subfield of $L$ containing the coefficients of
all the $p_{i}$, expressed as $L$-linear combinations of $B_{d}$. There are
only finitely many coefficients, and thus $K$ is a finite field extension of
$\f[p]$. Then, $X$ together with the $p_i$ yields a hsop of $R_K$.
\end{proof}

\begin{theorem}\label{thm:rel}
Let $R=\tau_nH^*(G;\f[p])$. Assume that there is a finite extension field $K$
of $\f[p]$ such that $R_K$ has a parameters $p_1,...,p_r$ in degrees $d_1,...,d_r$,
and denote $N = \sum_{i=1}^rd_i - \depth\left(H^*(U;\f[p])\right)$.  Assume
that $\alpha_n$ is surjective, and $n\ge N$.  Then, $\alpha_n$ is an
isomorphism, if and only if $P(R;t)\cdot\prod_{i=1}^r (1-t^{d_i})$ is a
polynomial of degree at most $N$.
\end{theorem}

\begin{proof}
Since $H^*(G;K) \cong K\otimes_{\f[p]}H^*(G;\f[p])$, we have
$R_K=\tau_nH^*(G;K)$, $P(R;t)=P(R_K;t)$ and
$P\left(H^*(G;\f[p]);t\right)=P\left(H^*(G;K);t\right)$. As a rational
function, $$\deg\left(P\left(H^*(G;\f[p]);t\right)\right)\le
\Reg\left(H^*(G;\f[p])\right)-\depth\left(H^*(G;\f[p])\right).$$
But by a result of Symonds on Castelnuovo--Mumford
regularity~\cite{Symonds:RegProof}, one has
$\Reg\left(H^*(G;\f[p])\right)=0$.
Moreover, $\depth(H^*(G;\f[p]))\ge\depth(H^*(U;\f[p]))$
by~\cite[Thm 2.1]{Benson:NYJM2}\footnote{Note that the statement originally is for a
  Sylow $p$-subgroup, but the proof only requires the index to be coprime
  to~$p$.}.
Hence, $$\deg\left(P\left(H^*(G;\f[p]);t\right)\right)\le -\depth(H^*(U;\f[p])).$$

Since $\alpha_n$ is surjective, the parameters $p_{i}$ of $R_K$ correspond to
parameters of $H^*(G;K)$ of the same degrees. Therefore and since
$P\left(H^*(G;K);t\right)$ is equal to $P\left(H^*(G;\f[p]);t\right)$, we
obtain that $p(t)=P\left(H^*(G;\f[p]);t\right)\cdot\prod_{i=1}^r (1-t^{d_i})$
is a polynomial. Since $P\left(H^*(G;\f[p]);t\right)$ is a rational function
of degree at most $-\depth\left(H^*(U;\f[p])\right)$, it follows that
$\deg(p(t))\le N$. This proves the ``only if'' part of the theorem.

Conversely, let $q(t) = P\left(R;t\right)\cdot\prod_{i=1}^r (1-t^{d_i})$. It
is easy to see that the degree-$\le N$ part of $q(t)$ is determined by
$\dim_{\f[p]}(R^{(d)})=\dim_K(R_K^{(d)})$ for $d=1,...,N$. But since $n\ge N$,
we have $\dim_{\f[p]}(R^{(d)})=\dim_{\f[p]}\left(H^{(d)}(G;\f[p])\right)$,
which in turn is determined by $p(t)$. Hence, if $q(t)$ is a polynomial of
degree at most $N$ then $q(t)=p(t)$ and thus
$P\left(R;t\right)=P\left(H^*(G;\f[p]);t\right)$. Since $\alpha_n$ is
surjective, this implies that $\alpha_n$ is an isomorphism.
\end{proof}

\section{Application of the criteria}\label{sec:expl}
\label{sec:example}

\subsection{On Lemma~\ref{lem:gen}}

When computing the $n$-th ring approximation of a cohomology ring with the
stable element method, then it is very handy to know surjectivity of
$\alpha_n$. Namely, identify $H^{(n)}(G,\f[p])$ with the stable subspace of
$H^{(n)}(U,\f[p])$ for an appropriate subgroup $U<G$. The stable subspace is
characterised by systems of linear equations in $H^{(n)}(U,\f[p])$. The number
of equations grows depends on the number of double cosets of $U$ in $G$, and
on the dimension of $H^{(n)}(U\cap U^c;\f[p])$, where $U^c$ denotes the
conjugate of $U$ under a double coset representative $c$.

Typically, if the last generator of a minimal ring presentation of
$H^{(n)}(G,\f[p])$ can be found in degree $d$, then the last relation can be
found in degree $2d$. And typically, the size of the equation systems grows
rapidly with the degree: In high degrees, the formulation and solution of the
stability conditions would often require a considerable amount of computation
time. This can be avoided with Lemma~\ref{lem:gen}. 

To compute the $n$-th ring approximation, one would usually start with the
power products of generators of $H^*(G,\f[p])$ in degrees less than $n$, which
allows to find new algebraic relations between these generators. Then, one
would compute $H^{(n)}(G,\f[p])$ by solving the stability conditions, would
compare it with the subspace spanned by power products, and may find further
generators of $H^*(G,\f[p])$.
But if $\alpha_n$ is surjective, then $H^{(n)}(G,\f[p])$ coincides with the
linear span of products of elements of smaller degree, which has been
determined to detect relations in degree $n$ anyway. Hence, for computing
$H^{(n)}(G,\f[p])$, one does not need to compute large systems of equations in
$H^{(n)}(U,\f[p])$.

As an example, we consider the cohomology ring of $\mathrm{SuzukiGroup}(8)$
with coefficients in $\f[2]$.  It turns out that this is isomorphic to the
cohomology ring of group number $179$ of order $448$ in the Small Groups
library~\cite{BeEiOBr:Millennium}, here denoted by $G$; this group is the
normaliser of the centre of a Sylow $2$--subgroup $S$ of
$\mathrm{SuzukiGroup}(8)$.  It turns out that six double cosets of $S$ in $G$
are enough to determine the stability conditions that describe $H^*(G;\f[2])$
as a graded sub-ring of $H^*(S;\f[2])$.

When computing a ring approximation of $H^*(G;\f[2])$ in increasing degrees,
the major part of the computation time is spent on the stability
conditions. One finds 99 generators out to degree $28$. It turns out that
$H^*(S;\f[2])$ is not a finitely generated module over $\tau_nH^*(G;\f[2])$,
for any $n<28$. But $H^*(S;\f[2])$ can be generated as a module over
$\tau_{28}H^*(G;\f[2])$ by elements of degree at most $29$. Hence,
$\gendeg_{28}(G,S) = 29$, and Lemma~\ref{lem:gen} asserts that a minimal ring
presentation of $H^*(G;\f[2])$ has no generators in degree $> 29$.

It turns out that indeed there are three generators in degree $29$. But in
degree $30$ and beyond, it is not needed to consider the stability
conditions. This saves enough resources to allow for a complete computation of
a minimal ring presentation of $H^*(G;\f[2])$, which is formed by $102$
generators of degree at most $29$, and $4790$ relations of maximal degree
$58$.

\subsection{Comparing the different criteria}

There are different criteria to prove completeness of an approximation of a
modular cohomology ring: The (modified) Benson
criterion, the Symonds criterion and the Hilbert--Poincar\'e criterion. 
It is not a-priori clear which criterion will be best or easiest to
use. Hence, given an approximation of a specific cohomology ring, it makes
sense to let a heuristics decide which criterion should be used to test
completeness. In this subsection, we discuss the reasoning behind such a
heuristics, namely the advantages and disadvantages of the different
criteria. Let $G$ be a finite group.

For the Benson and the Hilbert--Poincar\'e criteria, it helps to have a lower
bound for the depth of the cohomology ring $H^*(G,\f[p])$. In any case, the
depth is at least the rank of the centre of a Sylow $p$-subgroup of $G$, by
Duflot's theorem \cite[Thm 12.3.3]{CarlsonTownsley}. This bound is, of course,
easy to obtain.
If $G$ is not of prime power order and $H^*(G,\f[p])$ is computed by the
stable element method using a subgroup $U<G$, then the depth of $H^*(G,\f[p])$
is at least the depth of $H^*(U,\f[p])$. It can be difficult to compute the
depth of $H^*(U,\f[p])$. If it turns out to be too difficult, then one should
resort to Duflot's bound, although that might be weaker. In the following
paragraphs, let $D$ be a lower bound for the depth of $H^*(G,\f[p])$.

All criteria considered here involve degrees of parameters of modular
cohomology rings. Benson's criterion additionally needs parameters that are
filter-regular. How can a filter-regular hsop be found?

According to Benson, one can find elements that simultaneously restrict to
powers of the Dickson invariants in the cohomology rings of all maximal
$p$-elementary abelian subgroups. They form a filter-regular homogeneous
system of parameters. Since one has explicit formulae for the Dickson
invariants, the simultaneous lifts can be effectively constructed in a ring
approximation of sufficient degree. The Dickson invariants have degrees
$p^r-p^k$, where $r$ is the $p$-rank of $G$ and $k=0,...,r-1$.

In~\cite{128gps}, we proposed a slightly different construction that is always
available if $G$ is of prime power order: One can choose the generators of
$H^*(G,\f[p])$ so that they contain a sequence of elements that restrict to a
regular sequence in the cohomology of $\Omega_1(Z(Syl_p(G)))$, and there are
elements that simultaneously restrict to the Dickson invariants in the
cohomology rings of the complements of $\Omega_1(Z(Syl_p(G)))$ in the maximal
$p$-elementary abelian subgroups of $G$. The method is not guaranteed to work
for groups that are not of prime power order. But \emph{if} it works, then one
obtains a filter-regular hsop that is formed by $c$ generators together with
elements of degrees $p^{r-c}-p^k$, where $c$ is the rank of $Z(Syl_p(G))$, and
$k=0,...,r-c-1$. Moreover,~\cite{128gps} provides some ways to improve a given
f.-r.\@ hsop.

By now, we assume that we have constructed filter-regular parameters in degree
$d_1,...,d_r$. If one uses these parameters in the Hilbert--Poincar\'e
criterion, then it can only apply in degree $\ge d_1+...+d_r-D$. Using the
same parameters, Symonds' criterion can only apply in degree
$d_1+...+d_r-r+1$. The same holds for Benson's criterion, unless $D>1$: In
this case, Benson's criterion could already apply in degree $d_1+...+d_r-r$.
Hence, one keeps focusing on explicitly constructed filter-regular parameters,
it seems that Benson's criterion will always be best.

In~\cite{128gps}, we have shown how to modify Benson's by using an existence
proof for filter-regular parameters over a finite field extension. In many
cases, there exist a f.-r.\@ hsop in smaller degrees $\tilde d_1,...,\tilde
d_r$ (over a finite field extension) than the one that was explicitly
constructed. If $D>1$, the modified Benson criterion could potentially detect
completeness in degree $\tilde d_1+...+\tilde d_r-r$, and the
Hilbert--Poincar\'e criterion could apply in degree $\tilde d_1+...+\tilde
d_r-D$. Again, it seems that Benson's criterion will be best.

However, Symonds' criterion has the advantage that it can work with parameters
that do not form a filter-regular sequence and are not even algebraically
independent. We suggest two ways to construct such parameters.

The first approach starts with a filter-regular hsop. One could then mod out
all but one of the parameters, and try to find a smaller one by
enumeration. Of course, enumeration can be very expensive in general. But
Lemma~\ref{lem:pars} helps to restrict the search, so that enumeration often
becomes feasible. The result is a set of algebraically independent parameters,
but is not necessarily filter-regular.

The second approach starts with verifying that the cohomology rings of maximal
$p$-elementary abelian subgroups are finitely generated modules over the
restriction of $\tau_nH^*(G,\f[p])$. If this is the case, then $H^*(G,\f[p])$
is finite over the image of $\alpha_n$ by~\cite{Carlson:Tests}, and hence any
hsop in $\tau_nH^*(G,\f[p])$ maps to a hsop of $H^*(G,\f[p])$.
One tries to find a minimal subset $S$ of the generators of
$\tau_nH^*(G,\f[p])$ such that the quotient of $\tau_nH^*(G,\f[p])$ by the
sub-algebra spanned by $S$ is a finite-dimensional vector space over
$\f[p]$. Then, $S$ maps to a set of parameters of $H^*(G,\f[p])$ that, in
general, is algebraically dependent. Sometimes, the second method yields a
better degree sum for the parameters than the first method.

Let $e_1,...,e_n$ be the degrees of the parameters thus obtained; we have
$n\ge r$, and $n>r$ if and only if the parameters are algebraically
dependent. Symonds' criterion would potentially apply in degree
$e_1+...+e_n-n+1$, while the Hilbert-Poincar\'e criterion would potentially
apply in degree $e_1+...+e_n-D$. There are cases in which this is better than
the bound obtained with the modified Benson criterion.

Since $D\le r$ and often $r<n$, it would seem that typically Symonds'
criterion will be slightly better than the Hilbert--Poincar\'e
criterion. However, Lemma~\ref{lem:pars} may be able to prove the existence of
parameters in smaller degrees, over some extension field.
We found in practical computations: If Lemma~\ref{lem:pars} states the
existence of parameters of a certain degree over an extension field, then we
are usually to find such parameters \emph{without} a field extension, by
enumeration. But the point is that enumeration can be quite expensive. 
While the Hilbert--Poincar\'e criterion can exploit the mere existence of
parameters in small degrees, the Symonds criterion relies on their explicit
construction.

These considerations give rise to the following heuristics\footnote{The
  heuristics tries to minimize the effort needed to prove completeness, and is
  guaranteed to terminate in finite time} that is used in our optional Sage
package~\cite{SimonsProg}. Compute approximations of $H^*(G,\f[p])$ in
increasing degree, until $\tau_n H^*(G,\f[p])$ contains parameters for
$H^*(G,\f[p])$. Use Lemma~\ref{lem:gen} if $G$ is not of prime power order.
If there is a generator in degree $d$ that is not a regular element, then
compute the approximation at least out to degree $2d$. If this degree is
atteined, try to prove completeness (and terminate the computation, if
possible) as follows.

\begin{itemize}
\item Find parameters $S$ among the generators of $\tau_n H^*(G,\f[p])$. If
  their degrees are small enough for Symonds' criterion to apply, then try to
  prove completeness by that criterion.
\item Construct filter-regular parameters $F$ for $H^*(G,\f[p])$, by
  methods from~\cite{128gps}, say. Try to improve it by enumeration in small
  degrees (say, one and two), and obtain a potentially smaller but not
  necessarily filter-regular set $F'$ of parameters.
\item If the parameter degrees in $F'$ are small enough for Symonds' criterion
  to apply, then try to prove completeness by that criterion.
\item If $G$ is not of prime power order and Lemma~\ref{lem:pars} applied to
  either $S$ or $F'$ yields the existence of parameters in degrees small
  enough for the Hilbert--Poincar\'e criterion to apply, then try to prove
  completeness by that criterion.
\item Try to test completeness with the modified Benson criterion
  from~\cite{128gps}, but skip the test if the computation of the filter
  degree type turns out to be too difficult.
\end{itemize}

\bibliographystyle{abbrv}
\bibliography{references}

\end{document}